\documentclass[11pt]{article}

\usepackage{amssymb}

\usepackage{amsmath}
\allowdisplaybreaks

\newcommand{\rom}[1]{{\rm #1}}
\newtheorem{theorem}{Theorem}
\newtheorem{lemma}{Lemma}
\newtheorem{remark}{Remark}

\begin{document}

\begin{center}{\Large \bf
  On convergence of generators of equilibrium dynamics of hopping particles to generator of a birth-and-death process in continuum
}\vspace{2mm}

{\large \bf E. Lytvynov and P.T. Polara}\\
Department of Mathematics,  Swansea University, Singleton
Park, Swansea SA2 8PP, U.K.\end{center}

\begin{abstract}

We deal with two following classes of equilibrium stochastic dynamics of infinite particle systems in 
continuum: hopping particles (also called Kawasaki dynamics), i.e., a dynamics where each particle randomly hops over the space,
and birth-and-death process in continuum (or Glauber dynamics), i.e., a dynamics where there is no motion of particles, but 
rather particles die, or are born at random. We prove that a wide class of Glauber dynamics can be derived as a scaling limit 
of Kawasaki dynamics. More precisely, we prove the convergence of respective generators on a set of cylinder functions, in 
the  $L^2$-norm with respect to the invariant measure of the processes. The latter measure is supposed to be a Gibbs measure corresponding to a potential of pair interaction, in the low activity--high temperature regime. Our result 
generalizes that of [Finkelshtein~D.L. et al., to appear in Random Oper.\ Stochastic Equations], which was proved  for a special Glauber (Kawasaki, respectively) dynamics.\end{abstract}
\vspace{3mm}

 \vspace{1.5mm}

 \vspace{1.5mm}
%

\section{Introduction}

This paper deals with two classes of equilibrium stochastic dynamics of infinite particle systems in
continuum. Let $\Gamma$ denote the space of all locally finite subsets of $\mathbb R^d$. Such a space is called the 
configuration space (of an infinite particle system in continuum). Elements of $\Gamma$ are called configurations and each 
point of a configuration represents position of a particle.

One can naturally define a $\sigma$-algebra on $\Gamma$, and then a probability measure on $\Gamma$ represents a random 
system of particles. A probability measure on $\Gamma$ is often called a point process (see e.g.\ \cite{Kal}). Configuration 
 spaces and point processes are important tools of  classical statistical mechanics of continuous systems. A central class
of point processes which is studied there is the class of Gibbs measures. Typically one deals with Gibbs measures which 
correspond to a potential of pair interaction.

An equilibrium stochastic dynamics in continuum is a Markov process on $\Gamma$ which has a point process (typically a Gibbs measure) $\mu$ as
its invariant measure. 
One can distinguish three main classes of stochastic dynamics: 
\begin{itemize}

\item
diffusion processes, i.e., dynamics where each particle 
continuously moves in the space, see e.g.\ \cite{AKR4,Fritz,KLRDiffusions,MR98,Osa96,RS98,Yos96};

\item birth-and-death processes in continuum (Glauber dynamics), i.e., dynamics where 
there is no motion of particles, but rather particles disappear (die) or appear (are born) at random, see e.g.\ \cite{BCC,G1,HS,KL,KLR,KMZ,P,Wu};

\item  hopping particles (Kawasaki dynamics), i.e., dynamics where each 
particle randomly hops over the space \cite{KLR}. 
\end{itemize}

For a deep understanding of these dynamics, it is important to see how they are related to each other. For example, in the recent paper
\cite{KKL}, it was shown that a typical diffusion dynamics can be derived through a diffusive scaling limit of a 
corresponding Kawasaki dynamics.

In \cite{FKL}, it was proved that a special Glauber dynamics can be derived 
through a scaling limit of Kawasaki dynamics. Furthermore, \cite{FKL} conjectured that such a result holds, in fact, for a wide class of birth-and-death dynamics (dynamics of hopping particles, respectively), which are indexed by a parameter $s\in[0,1]$.
(Note that the result of \cite{FKL} corrrespods to the choice of parameter $s=0$.)

The  aim of this work is to show that the conjecture of \cite{FKL} is indeed true, at least for parameters $s\in[0,1/2]$. (In the case where $s\in (1/2,1]$, one needs to put additional, quite restrictive assumptions on the potential of pair interaction, and we will not treat this case in the present paper.)
 Thus, we show that the result of \cite{FKL} is not a  property of just one special Kawasaki (Glauber, respectively) 
dynamics, but rather represents a property which is common for many dynamics. 

More specifically, we fix a class of cylinder functions on $\Gamma$, and prove that on this class of functions, the 
corresponding generators converge in the $L^2(\Gamma,\mu)$-space. Here, $\mu$ is a Gibbs measure in the low activity--high temperature regime, $\mu$ being invariant measure for all the processes under consideration.   If one additionally knows that the class of cylinder 
functions is a core for the limiting generator, then our result implies weak convergence of finite-dimensional 
distributions of the corresponding processes. Unfortunately, apart from a very special case \cite{KL}, no result about a core for 
these generators is yet available.

The paper is organized as follows. In Section 2, we briefly discuss Gibbs measures in the low activity--high temperature regime, 
and the corresponding correlation and Ursell functions. In Section 3, we describe classes of  birth-and-death processes and  of dynamics of hopping particles. In Section~4, we formulate and prove the result about convergence of the generators.  

The authors acknowledge numerous useful discussions with Dmirti Finkelshtein and Yuri Kondratiev.

\section{Gibbs measures in the low activity-high temperature regime}

The configuration space over $\mathbb R^d$,  $d\in\mathbb N$,
is defined by 
\[
\Gamma:=\{\gamma\subset\mathbb R^d :\, |\gamma\cap\Lambda|<\infty \text{
for each compact } \Lambda\subset\mathbb R^d \},
\]
where $|\cdot|$ denotes the cardinality of a set. One can identify any $\gamma\in\Gamma$
with the positive Radon measure $\sum_{x\in\gamma}\varepsilon_x \in
\mathcal M(\mathbb R^d)$, where $\varepsilon_x$ is the Dirac measure with mass
at $x$, $\sum_{x\in\varnothing}\varepsilon_x:=$zero measure, and
$\mathcal M(\mathbb R^d)$ stands for the set of all positive Radon measures
on the Borel $\sigma$-algebra $\mathcal B(\mathbb R^d)$. The space $\Gamma$
can be endowed with the relative topology as a subset of
the space $\mathcal M(\mathbb R^d)$ with the vague topology, i.e., the
weakest topology on $\Gamma$ with respect to which all maps
\[
\Gamma\ni\gamma\mapsto\langle f,\gamma\rangle:=\int_{\mathbb R^d} f(x)\,\gamma(dx)=\sum_{x\in\gamma}f(x),\quad f\in C_0(\mathbb R^d),
\]
are continuous. Here, $C_0(\mathbb R^d)$ is
the space of all continuous real-valued functions on $\mathbb R^d$
with compact support. We will denote by $\mathcal B(\Gamma)$ the
Borel $\sigma$-algebra on $\Gamma$. 

Let $\mu$ be a probability measure on $(\Gamma,\mathcal B(\Gamma))$. Assume that, for each $n\in\mathbb N$, there
exists a non-negative, measurable symmetric function
$k_\mu^{(n)}$ on $(\mathbb R^d)^n$ such that, for any measurable
symmetric function $f^{(n)} :(\mathbb R^d)^n \to
[0,+\infty]$,
\[
\int_\Gamma \sum_{\{x_1,\dots,x_n\}\subset\gamma}
f^{(n)}
(x_1,\dots,x_n)\,\mu(d\gamma)
=\frac{1}{n!}\int_{(\mathbb R^d)^n}f^{(n)}
(x_1,\dots,x_n)k_\mu^{(n)}
(x_1,\dots,x_n)\,d x_1 \dotsm d x_n.
\]

The functions $k_\mu^{(n)}$ are called correlation functions
of the measure $\mu$. If there exists a constant $\xi>0$
such that
\begin{equation}\label{RB}
\forall(x_1,\dots,x_n)\in(\mathbb R^d)^n
:\quad k_\mu^{(n)} (x_1,\dots,x_n)\le\xi^n,
\end{equation}
then we say that the correlation functions $k_\mu^{(n)}$
satisfy the Ruelle bound.

The following lemma gives a characterization of the correlation functions in terms of the Laplace transform of a given point process, see e.g.\ \cite{KKL}

\begin{lemma}\label{n} Let $\mu$ be a probability measure on $(\Gamma,\mathcal B(\Gamma))$
which satisfies the Ruelle bound \eqref{RB}. Let $f:\mathbb{R}^d\to\mathbb{R}$ be a measurable function which is bounded outside a compact set 
$\Lambda \subset \mathbb R^d$ and such that $e^{f}-1 \in L^1(\mathbb R^d,dx)$. Then for $\mu$-a.a.\ $\gamma \in \Gamma$,
$\langle |f|,\gamma \rangle < \infty$ and
\[\int_\Gamma e^{\langle f,\gamma\rangle} \mu(d\gamma)=1+\sum_{n=1}^\infty \frac{1}{n!}\int_{(\mathbb{R}^d)^n} 
(e^{f(x_1)}-1)\dotsm (e^{f(x_n)}-1)k_\mu^{(n)}(x_1,\dots,x_n)\, dx_1\dotsm dx_n.\]
\end{lemma}

\begin{remark}\label{fff}\rom{Note that if $f:\mathbb{R}^d\to\mathbb{R}$ is bounded outside a compact set $\Lambda\subset\mathbb R^d$ and if, furthermore, $f$ is bounded from above on the whole $\mathbb R^d$, then the condition $e^{f}-1 \in L^1(\mathbb R^d,dx)$
is equivalent to $f\in L^1(\Lambda^c,dx)$.
}\end{remark}

Via a recursion formula, one can transform the correlation functions
$k_\mu^{(n)}$ into the Ursell functions $u_\mu^{(n)}$ and vice
versa, see e.g.\ \cite{Ru69}. Their relation is given by
\begin{equation}\label{urs_def}
k_\mu(\eta)=\sum u_\mu(\eta_1)\dotsm u_\mu(\eta_j),\quad
\eta\in\Gamma_0,\ \eta\ne\varnothing,
\end{equation}
where
\[
\Gamma_0:=\{\gamma\in\Gamma:  |\gamma|<\infty \},
\]
for any $\eta=\{x_1,\dots,x_n\}\in\Gamma_0$
\[
k_\mu(\eta):=k_\mu^{(n)}(x_1,\dots,x_n), \quad
u_\mu(\eta):=u_\mu^{(n)}(x_1,\dots,x_n) ,
\]
and the summation in (\ref{urs_def}) is over all partitions
of the set $\eta$ into nonempty mutually disjoint subsets
$\eta_1,\dots,\eta_j\subset\eta$ such that
$\eta_1\cup\dotsm\cup\eta_j=\eta$, $j\in\mathbb N$. 
Note that if the correlation functions $(k_\mu^{(n)})_{n=1}^\infty$ are translation invariant, i.e., for each $a\in\mathbb R^d$
\[ k_\mu^{(n)} (x_1,\dots,x_n)=k_\mu^{(n)}(x_1+a,\dots,x_n+a),\quad (x_1,\dots,x_n)\in(\mathbb R^d)^n,\]
then so are the Ursell functions $(u_\mu^{(n)})_{n=1}^\infty$.

A pair potential is a Borel-measurable function $\phi:\mathbb R^d \to\mathbb R\cup\{+\infty\}$
such that $\phi(-x)=\phi(x)\in\mathbb R$ for all $x\in\mathbb R^d\setminus\{0\}$.
For $\gamma\in\Gamma$ and $x\in\mathbb R^d\setminus\gamma$, we define a relative
energy of interaction between a particle at $x$ and the
configuration $\gamma$ as follows:
\[
E(x,\gamma):=\left\{
\begin{aligned}
&\sum_{y\in\gamma}\phi(x-y),&& \text{if } \sum_{y\in\gamma} |\phi(x-y)|
< +\infty ,\\
& +\infty, &&\text{otherwise.}
\end{aligned}
\right.
\]

A probability measure $\mu$ on $(\Gamma,\mathcal B(\Gamma))$ is
called a (grand canonical) Gibbs measure corresponding to
the pair potential $\phi$ and activity $z>0$ if it satisfies
the Georgii--Nguyen--Zessin identity (\cite[Theorem 2]{NZ}):
\begin{equation}\label{GNZ}
\int_\Gamma\mu(d\gamma)\int_{\mathbb R^d}\gamma(dx) F(\gamma,x)
=\int_\Gamma \mu(d\gamma)\int_{\mathbb R^d} z\,dx \exp[-E(x,\gamma)]
F(\gamma\cup x,x)
\end{equation}
for any measurable function $F:\Gamma\times\mathbb R^d\to[0;+\infty]$. Here and below, for simplicity of notations, we just write $x$ instead of $\{x\}$.
We denote the set of all such measures $\mu$ by $\mathcal G (z,\phi)$.

As a straightforward corollary of the Georgii--Nguyen--Zessin
identity (\ref{GNZ}), we get the following equality:
\begin{align}\notag
&\int_\Gamma \mu(d\gamma)\int_{\mathbb R^d}\gamma(dx_1)\int_{\mathbb R^d}\gamma(dx_2)
F(\gamma,x_1,x_2)\\
&=\int_\Gamma \mu(d\gamma)\int_{\mathbb R^d}z\, dx_1\int_{\mathbb R^d} z\, dx_2\exp\left[
-E(x_1,\gamma)-E(x_2,\gamma)-\phi(x_1-x_2)
\right]\notag\\
&\qquad \times F(\gamma\cup \{x_1, x_2\},x_1,x_2)\notag\\
&\quad + \int_\Gamma \mu(d\gamma)\int_{\mathbb R^d} z\,dx \exp\left[
-E(x,\gamma)\right] F(\gamma\cup x,x,x)\label{double_GNZ}
\end{align}
for any measurable function
$F:\Gamma\times\mathbb R^d\times\mathbb R^d\to[0,+\infty]$.

Let us formulate  conditions on the pair potential $\phi$.

{\bf (S) (Stability)}
 There exists $B\ge0$ such that, for any $\gamma\in\Gamma_0$, 
\[
\sum_{\{x,y\}\subset\gamma}\phi(x-y)\ge -B|\gamma|.
\]

In particular, condition (S) implies that $\phi(x)\geq -2B$,
$x\in\mathbb R^d$.

{\bf (LA-HT) (Low activity-high temperature regime)} We have:
\[
\int_{\mathbb R^d}|e^{-\phi(x)}-1|z\,dx<(2e^{1+2B})^{-1},
\]
where $B$ is as in (S).

The following classical theorem is due to Ruelle
\cite{wewewe,Ru69}.

\begin{theorem}\label{cfdrer} Assume that \rom{(S)} and \rom{(LAHT)} are satisfied. Then there exists $\mu\in\mathcal G(z,\phi)$ which has the following properties\rom:

\rom{a)} $\mu$ has correlation functions $(k_\mu^{(n)})_{n=1}^\infty$\rom, which are translation invariant
and satisfy the Ruelle bound \eqref{RB}\rom;

\rom{b)} For each $n \ge 2$, we have $u_\mu^{(n)}(0,\cdot,\cdot,\dots,\cdot)\in
L^1(\mathbb R^{d(n-1)},dx_1 \dotsm dx_{n-1})$, where \linebreak $u_\mu^{(n)}(0,\cdot,\cdot,\dots,\cdot)$
is considered as a function of $n-1$ variables.
\end{theorem}

In what follows, we will assume that (S) and (LAHT) are satisfied, and we will 
keep the measure $\mu$ from Theorem~\ref{cfdrer} fixed. 

\section{Equilibrium birth-and-death (Glauber) dynamics  and hopping particles' \newline
(Kawasaki) dynamics}

In what follows,  we will additionally assume that $\phi$ is bounded outside some ball in $\mathbb R^d$. Note that then (see e.g.\ \cite{KLR})
$$E(x,\gamma)=\sum_{y\in\gamma}\phi(x-y),$$ for $dx\,\mu(d\gamma)$-a.a. $x\in\mathbb R^d$ and $\gamma\in\Gamma$ and
$$E(X,\gamma\setminus x)=\sum_{y\in\gamma\setminus x}\phi(x-y),$$ for $\mu$-a.a. $\gamma\in\Gamma$ and 
all $x\in\gamma$.
 
We fix a parameter $s \in [0,1/2]$. We introduce the set $\mathcal FC_b(C_0(\mathbb R^d),\Gamma)$ of all functions of the form 
$$\Gamma \ni \gamma \mapsto F(\gamma)=g(\langle f_1,\gamma \rangle,\dots,\langle f_N,\gamma \rangle),$$
where $N\in\mathbb N$, $f_1,\dots,f_N \in C_0(\mathbb R^d)$, and $g\in C_b(\mathbb R^N)$, where $C_b(\mathbb R^N)$
denotes the set of all continuous bounded functions on $\mathbb R^N$. For each function $F:\Gamma \to \mathbb R$,
$\gamma\in\Gamma$, and $x,y\in\mathbb R^d$, we denote
$$(D_{x}^{-}F)(\gamma):=F(\gamma \setminus x)-F(\gamma),$$
$$(D_{xy}^{-+}F)(\gamma):=F(\gamma \setminus x \cup y)-F(\gamma).$$

We fix a bounded function $a:\mathbb R^d \to [0,+\infty)$ such that $a(-x)=a(x)$, $x\in\mathbb R^d$, and $a\in L^1(\mathbb R^d,dx)$.
We define bilinear forms
\begin{align*}
\mathcal E_G(F,G)&=\int_\Gamma\mu (d\gamma)\int_{\mathbb R^d} \gamma(dx) 
\exp[s E(x,\gamma\setminus x)](D_{x}^{-}F)(\gamma)(D_{x}^{-}G)(\gamma),
\end{align*}
\begin{align*}
\mathcal E_K(F,G)&=\frac{1}{2}\int_\Gamma\mu (d\gamma)\int_{\mathbb R^d} \gamma(dx)\int_{\mathbb R^d}dy \ a(x-y)\\  
&\quad \times\exp[s E(x,\gamma\setminus x)-(1-s)E(y,\gamma \setminus x)](D_{xy}^{-+}F)(\gamma)(D_{xy}^{-+}G)(\gamma),\end{align*}
where $F,G \in\mathcal FC_b(C_0(\mathbb R^d),\Gamma)$. As we will see below, $\mathcal E_G$
corresponds to a Glauber dynamics and $\mathcal E_K$   corresponds to a Kawasaki dynamics.

The next theorem follows from \cite{KLR}.

\begin{theorem}\label{Hunt_thm}
\rom{i)} The bilinear forms $(\mathcal E_G,\mathcal FC_b(C_0(\mathbb R^d),\Gamma))$ and $(\mathcal E_K,\mathcal FC_b(C_0(\mathbb R^d),\Gamma))$ are closable on $L^2(\Gamma,\mu)$
and their closures are denoted by $(\mathcal E_G,D(\mathcal E_G))$ and  $(\mathcal E_K,D(\mathcal E_K))$, respectively. 

\rom{ii)} Denote by $(H_G,D(H_G))$ and  $(H_K,D(H_K))$ the generators of $(\mathcal E_G,D(\mathcal E_G))$ and \linebreak $(\mathcal E_K,D(\mathcal E_K))$, respectively. Then $\mathcal FC_b(C_0(\mathbb R^d),\Gamma)\subset D(H_G)\cap D(H_K)$, and for any $F\in \mathcal FC_b(C_0(\mathbb R^d),\Gamma)$, 
\begin{align}
(H_GF)(\gamma)&=-\int_{\mathbb R^d}\gamma(dx) \exp[sE(x,\gamma\setminus x)](D_{x}^{-}F)(\gamma)\notag\\&\qquad\text{}
-\int_{\mathbb R^d}z\, dx \exp[(s-1)E(x,\gamma)](D_{x}^{+}F)(\gamma),\label{fguftftf}\\
(H_KF)(\gamma)&= -\int_{\mathbb R^d} \gamma(dx)\int_{\mathbb R^d}dy \, a(x-y) 
 \exp[sE(x,\gamma \setminus x)+(s-1) E(y,\gamma\setminus x)](D_{xy}^{-+}F)(\gamma).\label{C} 
\end{align}

\rom{iii)} Let $\sharp:=G,K$. There exists a conservative
Hunt process
\[
\mathbf {M}^\sharp=(\mathbf{\Omega}^\sharp, \, \mathbf{F}^\sharp, \
(\mathbf {F}^\sharp_t)_{t\geq 0}, \, (\mathbf{\Theta}^\sharp_t)_{t\geq 0}, \, (\mathbf {X}^\sharp(t)_{t\geq 0}, \, (\mathbf{P}^\sharp_\gamma)_{\gamma\in\Gamma})
\]
on $\Gamma$ \rom(see e.g.  \rom{\cite[p.~92]{MaRo})} which is properly associated
with $(\mathcal E_\sharp,D(\mathcal E_\sharp))$, i.e., for all \rom($\mu$-versions
of\,\rom) $F\in L^2(\Gamma,\mu)$ and all $t>0$ the function
\[
\Gamma\ni\gamma\mapsto(p_t^\sharp F)(\gamma)
:=\int_{\mathbf\Omega^\sharp} F(\mathbf {X}^\sharp(t))d\mathbf {P}^\sharp_\gamma
\]
is an $\mathcal E_\sharp$-quasi-continuous version of
$\exp[-tH_\sharp]F$. $\mathbf{M}^\sharp$ is up to
$\mu$-equivalence unique \rom(cf.\ \rom{\cite[Chap.~IV, Sect.~6]{MaRo})}.
In particular, $\mathbf {M}^\sharp$ has $\mu$ as  invariant
measure.
\end{theorem}

\begin{remark}\rom{
In Theorem~\ref{Hunt_thm}, $\mathbf{M}^\sharp$ can be taken canonical,
i.e., $\mathbf{\Omega}^\sharp$ is the set
$D([0,+\infty),\Gamma)$ of all {\it c\'adl\'ag\/} functions
$\omega:\left[0,+\infty\right)\to\Gamma$ (i.e., $\omega$
is right continuous on $\left[0,+\infty\right)$ and has
left limits on $(0,+\infty)$);
$\mathbf{X}^\sharp(t)(\omega)=\omega(t)$, $t\geq0$,
$\omega\in\mathbf {\Omega}^\sharp$; $(\mathbf {F}^\sharp_t)_{t\geq0}$,
together with $\mathbf {F}^\sharp$, is the corresponding minimum
completed admissible family (cf.\ \cite[Section~4.1]{Fu80}); $\mathbf {\Theta}^\sharp_t$, $t\geq0$, are the corresponding
natural time shifts.
}\end{remark}

It follows from \eqref{fguftftf} that $H_G$ is (at least heuristically) the generator of a birth-and-death process, in which the factor $\exp[sE(x,\gamma\setminus x)]$ describes the rate at which particle $x$ of the configuration $\gamma$ dies, whereas the factor $\exp[(s-1)E(x,\gamma)]$ describes the rate at which, given a configuration $\gamma$, a new particle is born at $x$. We see that particles tend to die in high energy regions, i.e., if $E(x,\gamma\setminus x)$
is high, and they tend to be born al low energy regions, i.e., if $E(x,\gamma)$ is low.

Next, by \eqref{C}, $H_K$ is (again at least heuristically) the generator of a hopping particle dynamics, in which the factor
\[\exp[sE(x,\gamma\setminus x)+(s-1)E(y,\gamma\setminus x)]\] describes the rate at which a particle $x$ of configuration $\gamma$ hops to $y$. 
We see that this rate is high if the relative energy of interaction between $x$ and the rest of the configuration, $\gamma\setminus x$, is high, whereas the relative energy of interaction between $y$ and $\gamma\setminus x$ is low, i.e., particles tend to hop from high energy regions to low energy regions. 

\section{Scaling limit}

In this section, we will show that the birth-and-death dynamics considered in Section~3 can be treated as a limiting dynamics of hopping particles.
In other words, we will perform a scaling of Kawasaki dynamics which will lead to the Glauber dynamics. 
 We will
only discuss this convergence at the level of convergence of the generators on an appropriate set of cylinder 
functions. In fact, such a convergence implies weak convergence of finite-dimensional distributions of corresponding 
equilibrium processes if additionally the set of test functions forms a core for the limiting generator. However, in the general case, 
no core of  this generator is yet known and this is an open, important problem, which we hope to return to in our future 
research. Our results will hold for all $s\in[0,1/2]$ (see Section~3). They will generalize Theorem 4.1 in 
\cite{FKL}, which was proved in the special case $s=0$, and confirm the conjecture formulated in Section~6 of that paper.

So, let us consider the following scaling of the Kawasaki dynamics (for a fixed $s\in[0,1/2]$). Recall that, for each bounded function
$a:\mathbb R^d \to \mathbb R$  such that   $a(x)\ge0$, $a\in L^1(\mathbb R^d,dx)$, and 
$a(-x)=a(x)$ for all $x\in\mathbb R^d$, we have constructed the corresponding generator of the Kawasaki dynamics. We now fix 
an arbitrary $\varepsilon >0$ and define a function $a_\varepsilon:\mathbb R^d \to \mathbb R$ by 
$$a_\varepsilon(x)=\varepsilon^d a(\varepsilon x),\quad x\in \mathbb R^d.$$ Note that
$$\int_{\mathbb R^d} a_\varepsilon(x)\,dx=\int_{\mathbb R^d} a(x)\,dx.$$
By the properties of the function $a$, we evidently have that the function $a_\varepsilon$ is also bounded, satisfies $a_\varepsilon(x)\ge0$,
for all $x\in\mathbb R^d$, $a_\varepsilon\in L^1(\mathbb R^d,dx)$, and 
$a_\varepsilon(-x)=a_\varepsilon(x)$ for all $x\in\mathbb R^d$. Hence, we can construct the Kawasaki 
generator which corresponds to the function $a_\varepsilon$. It is convineant for us to denote this generator by 
$(H_\varepsilon, D(H_\varepsilon))$. We will also denote the generator of the Glauber dynamics by 
$(H_0, D(H_0))$. We first need the following lemma, whose proof is completely analogous to the proof of Lemma 4.1
in \cite{FKL}.

\begin{lemma}
For any $\varepsilon \ge 0$ and any $\varphi \in C_0(\mathbb R^d)$, the function $F(\gamma):=e^{\langle\varphi,\gamma\rangle}
$ belongs to $ D(H_\varepsilon)$ and the action of $H_\varepsilon$ on $F$ is given by the right hand side of 
formula \eqref{C} for $\varepsilon>0$ \rom(with $a$ replaced by $a_\varepsilon$\rom), respectively by the right hand side of \eqref{fguftftf} for $\varepsilon=0$.  
\end{lemma}

\begin{remark}\rom{For each $\varepsilon\ge0$, denote by $(\mathcal E_\varepsilon,D(\mathcal E_\varepsilon))$ the Dirichlet form with the generator $(H_\varepsilon,D(H_\varepsilon))$.
It can be easily proved that the set $\big\{\exp[\langle\varphi,\cdot\rangle]:\varphi\in C_0(\mathbb R^d)\big\}$ 
is dense in the Hilbert space $D(\mathcal E_\varepsilon)$ equipped with inner product $(F,G)_{D(\mathcal E_\epsilon)}:=\mathcal E(F,G)+(F,G)_{L^2(\Gamma,\mu)}$. 
}\end{remark}

We have
\begin{align*}
(D_{xy}^{-+}F)(\gamma)&=F(\gamma \setminus x \cup y)-F(\gamma)\\
&=-F(\gamma)+F(\gamma \setminus x)-F(\gamma \setminus x)+F(\gamma \setminus x \cup y)\\
&=(D_{x}^{-}F)(\gamma)+(D_{y}^{+}F)(\gamma \setminus x).
\end{align*}
So, we may rewrite the action of $H_\varepsilon$ for $\varepsilon > 0$ as follows:
$$H_\varepsilon :=H_\varepsilon^{+} + H_\varepsilon^{-},$$ where
\begin{align*}
(H_\varepsilon^-F)(\gamma)&= -\int_{\mathbb R^d} \gamma(dx)(D_{x}^{-}F)(\gamma)\int_{\mathbb R^d}dy \ a_\varepsilon(x-y) 
\exp[sE(x,\gamma \setminus x)+(s-1) E(y,\gamma\setminus x)] 
\end{align*}
and
\begin{align*}
(H_\varepsilon^+F)(\gamma)&= -\int_{\mathbb R^d} \gamma(dx)\int_{\mathbb R^d}dy \ a_\varepsilon(x-y) 
\exp[sE(x,\gamma \setminus x)+(s-1) E(y,\gamma\setminus x)](D_{y}^{+}F)(\gamma \setminus x). 
\end{align*}
We can also rewrite $$H_0: =H_0^{+} + H_0^{-},$$where
\begin{align*}
(H_0^-F)(\gamma)&=-\int_{\mathbb R^d}\gamma(dx) \exp[sE(x,\gamma\setminus x)](D_{x}^{-}F)(\gamma)
\end{align*}
and
$$(H_0^+F)(\gamma)=-\int_{\mathbb R^d}z\, dx \exp[-(1-s)E(x,\gamma)](D_{x}^{+}F)(\gamma).$$

\begin{theorem}
Let $s\in[0,1/2]$ be fixed. Assume that the pair potential $\phi$ and activity $z>0$ satisfy conditions 
\rom{(S)} and \rom{(LA-HT)}. Assume that $\phi$ is bounded outside some compact set in $\mathbb R^d$. Assume also that \begin{equation}\label{12345}
\phi(x) \to 0 \text{ as }|x| \to \infty.
\end{equation}
 Let $\mu$ be the Gibbs measure from 
$\mathcal G(z,\phi)$ as in Theorem \rom{\ref{cfdrer}}. Assume that the function $a$ is chosen so that
\begin{align} 
\int_{\mathbb R^d}a(x)dx=\bigg(\int_\Gamma \exp\bigg[(s-1) \sum_{u\in\gamma}\phi(u)\bigg]\mu(d\gamma)\bigg)^{-1}.\label{2}
\end{align}
Then, for any $\varphi \in C_0(\mathbb R^d)$, 
$$H_\varepsilon^{\pm}e^{\langle \varphi,\cdot\rangle} \to H_0^{\pm}e^{\langle \varphi,\cdot\rangle} \text{ in } 
L^2(\Gamma,\mu) \text{ as } 
\varepsilon \to 0,$$ so that
$$H_\varepsilon e^{\langle \varphi,\cdot\rangle} \to H_0 e^{\langle \varphi,\cdot\rangle} \text{ in } 
L^2(\Gamma,\mu) \text{ as } \varepsilon \to 0.$$
\end{theorem}

\begin{remark}\rom{
In fact, condition \eqref{12345} can be omitted, and instead one can use the fact that $\phi$ is an integrable function outside a compact set in $\mathbb R^d$ (compare with \cite{FKL}). However, in any reasonable application, the potential $\phi$ does satisfy condition \eqref{12345}. 
}\end{remark}

\begin{remark}\rom{
Note that the integral on the right hand side of \eqref{2} is well defined and finite due to Lemma~\ref{n}, see also Remark~\ref{fff}
}\end{remark}

\noindent {\it{Proof.}} 
We first need the following lemma, which generalizes Lemma 4.2 in \cite{FKL}.
\begin{lemma}\label{T}
Let a function $\psi: \mathbb R^d \to \mathbb R$ be such that $e^\psi-1$ is bounded and integrable. Suppose that $A\ge0$,
$B\ge0$, $x_1,x_2,y_1,y_2\in \mathbb R^d$ and $x_1 \neq y_1$. Then 
\begin{multline*}\int_\Gamma \exp\bigg[-A E\bigg(\frac{x_1}{\varepsilon}+x_2,\gamma\bigg)-B E\bigg(\frac{y_1}{\varepsilon}+y_2,\gamma\bigg)+
\langle \psi,\gamma\rangle\bigg]\mu(d\gamma)\\
\to \int_\Gamma \exp\bigg[-A\sum_{u\in\gamma}\phi(u)\bigg]\mu(d\gamma) \int_\Gamma \exp\bigg[-B\sum_{u\in\gamma}\phi(u)\bigg]\mu(d\gamma) 
\int_\Gamma \exp[\langle \psi,\gamma \rangle]\mu(d\gamma)\end{multline*}
as $\varepsilon \to 0$.
\end{lemma}

\noindent{\it Proof.} By Lemma \ref{n},
\begin{align}
&\int_\Gamma \exp\bigg[-A E((x_1/\varepsilon)
+x_2,\gamma)-B E((y_1/\varepsilon)+y_2,\gamma)+
\langle \psi,\gamma\rangle\bigg]\mu(d\gamma)\notag \\
&\qquad=1+ \sum_{n=1}^\infty \frac{1}{n!}\int_{(\mathbb R^d)^n}(\exp[-A\phi(\cdot -x(\varepsilon))-B\phi(\cdot -y(\varepsilon))
+ \psi(\cdot)]-1)^\otimes(u_1,\dots,u_n)\notag \\
&\qquad\quad\times k_\mu^{(n)}(u_1,\dots,u_n)du_1\dotsm du_n,\label{M}
\end{align}
where $x(\varepsilon):=(x_1/\varepsilon)+x_2$, $y(\varepsilon):=(y_1/\varepsilon)+y_2$.

Using the Ruelle bound, semi-boundedness of $\phi$ from below and the integrability of $\phi$ outside a compact set,
 we conclude  from the dominated convergence theorem that, in order to find the limit of  the 
right hand side of \eqref{M} as $\varepsilon \to 0$, it suffices to find the limit of each term
\begin{align}
C_\varepsilon^{(n)}:&=\int_{(\mathbb R^d)^n}(\exp[-A\phi(\cdot -x(\varepsilon))-B\phi(\cdot -y(\varepsilon))
+ \psi(\cdot)]-1)^{\otimes n}(u_1,\dots,u_n)
\notag \\
&\quad\times 
k_\mu^{(n)}(u_1,\dots,u_n)\,du_1\dotsm du_n \notag \\
&=\sum_{n_1+n_2+n_3=n} {n \choose n_1\, n_2\, n_3}\int_{(\mathbb R^d)^n}(f_{1,\varepsilon}^{\otimes n_1}
\otimes f_{2,\varepsilon}^{\otimes n_2}\otimes f_{3,\varepsilon}^{\otimes n_3})(u_1,\dots,u_n) \notag\\
&\quad\times k_\mu^{(n)}(u_1,\dots,u_n)\,du_1\dotsm du_n,\label{N}
\end{align}
where
\begin{align*} 
&f_{1,\varepsilon}(u):=(\exp[-\psi(u)]-1)\exp[-A\phi(u -x(\varepsilon))-B\phi(u -y(\varepsilon))],\\
&f_{2,\varepsilon}(u):=(\exp[-A\phi(u -x(\varepsilon))]-1) \exp[-B\phi(u -y(\varepsilon))],\\ 
&f_{3,\varepsilon}(u):=\exp[-B\phi(u -y(\varepsilon))]-1,\quad u\in\mathbb R.
\end{align*}

Using definition of  Ursell functions,
we see that
\begin{align*}
&\int_{(\mathbb R^d)^n}(f_{1,\varepsilon}^{\otimes n_1}
\otimes f_{2,\varepsilon}^{\otimes n_2}\otimes f_{3,\varepsilon}^{\otimes n_3})(u_1,\dots,u_n)
 k_\mu^{(n)}(u_1,\dots,u_n)\,du_1\dotsm du_n \\
&=\sum\int_{(\mathbb R^d)^n}(f_{1,\varepsilon}^{\otimes n_1}
\otimes f_{2,\varepsilon}^{\otimes n_2}\otimes f_{3,\varepsilon}^{\otimes n_3})(u_1,\dots,u_n) 
u_\mu(\eta_1)\dotsm u_\mu(\eta_j)\,du_1\dotsm du_n,
\end{align*}
where the summation is over all partitions $\{\eta_1,\dots,\eta_j\}$ of $\eta=\{u_1,\dots,u_n\}$. We now have to distinguish the three following cases.

Case 1: Each element $\eta_i$ of the partition is either a subset of $\{u_1,\dots,u_{n_1}\}$, or a subset of 
$\{u_{n_1+1},\dots,u_{n_1+n_2}\}$, or a subset of $\{u_{n_1+n_2+1},\dots,u_n\}$.
Set $$u_i'=u_i-x(\varepsilon),\quad i=n_1+1,\dots,n_2,$$
$$u_i'=u_i-y(\varepsilon),\quad i=n_2+1,\dots,n.$$
Then using the translation invariance of the Ursell functions, we get that the corresponding term is equal to
\begin{align}   
\int_{(\mathbb R^d)^n}(f_{1,\varepsilon}^{\otimes n_1}
\otimes g_{2,\varepsilon}^{\otimes n_2}\otimes g_{3,\varepsilon}^{\otimes n_3})(u_1,\dots,u_n) 
u_\mu(\eta_1)\dotsm u_\mu(\eta_j)\,du_1\dotsm du_n,\label{O}
\end{align}
where 
\begin{align*} 
&g_{2,\varepsilon}(u):=(\exp[-A\phi(u)]-1)\exp[-B\phi(u +((x_1-y_1)/\varepsilon)+x_2-y_2)],\\
&g_{3,\varepsilon}(u):=\exp[-B\phi(u)]-1,\quad u\in\mathbb R.
\end{align*}
Note that $x_1-y_1\not=0$ and so for any fixed $u$ (and $x_2,y_2$), we have
$$|u+((x_1-y_1)/\varepsilon)+x_2-y_2| \to +\infty \quad\text{as } \varepsilon \to 0.$$
By \eqref{12345} and the dominated convergence theorem, we therefore have that \eqref{O}
 converges to
\begin{align*} 
&\int_{(\mathbb R^d)^n}(\exp[-\psi(\cdot)]-1)^{\otimes n_1}\otimes (\exp[-A\phi(\cdot)]-1)^{\otimes n_2}\\
&\qquad \otimes (\exp[-B\phi(\cdot)]-1)^{\otimes n_3}(u_1,\dots,u_n)
u_\mu(\eta_1)\dotsm u_\mu(\eta_j)\,du_1\dotsm du_n.
\end{align*}

Case 2: There is an element of the partition which has non-empty intersections with both sets $\{u_1,\dots,u_{n_1}\}$
and $\{u_{n_1+1},\dots,u_n\}$.

Using Theorem \ref{cfdrer}, we have that, for each $n\in\mathbb N$,
$$U_\mu^{(n+1)}\in L^1((\mathbb R^d)^n,dx_1 \dotsm dx_n),$$
where
$$U_\mu^{(n+1)}(x_1,\dots,x_n):=u_\mu^{(n+1)}(x_1,\dots,x_n,0),\quad (x_1,\dots,x_n)\in(\mathbb R^d)^n.$$
Consider the integral
\begin{align*}
&\int_{(\mathbb R^d)^k}(\exp[-\psi(u_1)]-1)u_\mu^{(k)}(u_1,\dots,u_k)\,du_1 \dotsm du_k \\
&\qquad=\int_{(\mathbb R^d)^k}(\exp[-\psi(u_1)]-1)u_\mu^{(k)}(0,u_2-u_1,u_3-u_1,\dots,u_k-u_1)\,du_1 \dotsm du_k,
\end{align*}
where we used translation invariance of Ursell function.
By changing variables $u_1'=u_1,\,u_2'=u_2-u_1,\,\dots,\,u_k'=u_k-u_1$, we continue as follows:
\begin{align*}
&=\int_{(\mathbb R^d)^k}(\exp[-\psi(u_1')]-1)u_\mu^{(k)}(0,u_2',u_3',\dots,u_k')\,du_1' \dotsm du_k'\\
&=\int_{\mathbb R^d}(e^{-\psi(u_1)}-1)du_1 \times 
\int_{(\mathbb R^d)^{k-1}}U_\mu^{(k)}(u_2,u_3,\dots,u_k)du_2 \dotsm du_k.
\end{align*}
Note also that $$|x(\varepsilon)| \to +\infty \text{ and }|y(\varepsilon)| \to +\infty \text{ as }\varepsilon \to 0,$$
and hence, for each fixed $u\in\mathbb R^d$
$$\exp[-A\phi(u-x(\varepsilon))]-1 \to 0,\quad \exp[-B\phi(u-y(\varepsilon))]-1 \to 0$$
as $\varepsilon \to 0$. From here, using the dominated convergence theorem, we conclude that
\begin{align*} 
&\int_{(\mathbb R^d)^n}(f_{1,\varepsilon}^{\otimes n_1}
\otimes f_{2,\varepsilon}^{\otimes n_2}\otimes f_{3,\varepsilon}^{\otimes n_3})(u_1,\dots,u_n) \\
&\quad\times u_\mu(\eta_1)\dotsm u_\mu(\eta_j)\,du_1 \dotsm du_n \to 0 
\end{align*}
as $\varepsilon \to 0$.

Case 3: Case 2 is not satisfied, but there is an element $\eta_l$ of the partition which has non-empty 
intersections with both sets
$\{u_{n_1+1},\dots,u_{n_1+n_2}\}$, and $\{u_{n_1+n_2+1},\dots,u_n\}$.

Shift all the variables entering $\eta_l$ by $x(\varepsilon)$. Now, since $\exp[-A\phi]-1 \in L^1(\mathbb R^d ,dx)$, 
analogously to case 2, the term converges to zero 
as $\varepsilon \to 0$.

Thus, again using the definition of the Ursell functions, we get, for each $n\in\mathbb N$,
\begin{align*}
C_\varepsilon^{(n)} & \to \sum_{n_1+n_2+n_3=n} {n \choose n_1\, n_2\, n_3}\\
&\quad\times
\int_{(\mathbb R^d)^{n_1}} 
(\exp[\psi( \cdot )]-1)^{\otimes n_1}(u_1,\dots,u_{n_1}) k_{\mu}^{(n_1)}
(u_1,\dots,u_{n_1})
\,du_1 \dotsm du_{n_1}\\
&\quad\times\int_{(\mathbb R^d)^{n_2}} (\exp[-A\phi(\cdot)]-1)^{\otimes n_2}(u_{n_1+1},\dots,u_{n_1+n_2})k_\mu^{(n_2)}
(u_{n_1+1},\dots,u_{n_1+n_2})\\
&\qquad\qquad\quad\times du_{n_1+1} \dotsm du_{n_1+n_2}\\
&\quad\times\int_{(\mathbb R^d)^{n_3}} (\exp[-B\phi(\cdot)]-1)^{\otimes n_3}(u_{n_1+n_2+1},\dots,u_n )k_\mu^{(n_3)}
(u_{n_1+n_2+1},\dots,u_n)\\
&\qquad\qquad\quad\times
 du_{n_1+n_2+1} \dotsm du_n.
\end{align*}
Therefore, the right hand side of \eqref{M} converges to
\begin{align*}
&\bigg(1+ \sum_{n=1}^\infty \frac{1}{n!}\int_{(\mathbb R^d)^n}(\exp[-A\phi(\cdot)]-1)^{\otimes n}(u_1,\dots,u_n)
k_\mu^{(n)}(u_1,\dots,u_n)\,du_1\dotsm du_n \bigg)\\
&\times\bigg( 1+ \sum_{n=1}^\infty \frac{1}{n!}\int_{(\mathbb R^d)^n}(\exp[-B\phi(\cdot)]-1)^{\otimes n}(u_1,\dots,u_n)
k_\mu^{(n)}(u_1,\dots,u_n)\,du_1\dotsm du_n \bigg)\\
&\times \bigg(1+ \sum_{n=1}^\infty \frac{1}{n!}\int_{(\mathbb R^d)^n}(\exp[\psi(\cdot)]-1)^{\otimes n}(u_1,\dots,u_n)
 k_\mu^{(n)}(u_1,\dots,u_n)\,du_1\dotsm du_n\bigg)\\
&\qquad=\int_\Gamma \exp\bigg[-A\sum_{u\in\gamma}\phi(u)\bigg]\mu(d\gamma) \int_\Gamma \exp\bigg[-B\sum_{u\in\gamma}\phi(u)\bigg]\mu(d\gamma) 
\int_\Gamma \exp[\langle \psi,\gamma \rangle]\mu(d\gamma).
\end{align*}
as $\varepsilon \to 0,$
 which proves the lemma.\quad $\square$
 
Now we are in position to prove the theorem. We fix any $\varphi \in C_0(\mathbb R^d)$ and denote 
$F(\gamma):=e^{\langle \varphi,\gamma\rangle}$. It suffices to prove that 
\begin{equation}\label{D}
\int_\Gamma (H_\varepsilon^{\pm}F)^2(\gamma)\mu(d\gamma) \to \int_\Gamma (H_0^{\pm}F)^2(\gamma)\mu(d\gamma)
\text{ as } \varepsilon \to 0,
\end{equation}
\begin{equation}\label{S}
\int_\Gamma(H_\varepsilon^{\pm}F)(\gamma)(H_0^{\pm}F)(\gamma)\mu(d\gamma) \to \int_\Gamma (H_0^{\pm}F)^2(\gamma)\mu(d\gamma)
\text{ as } \varepsilon \to 0.
\end{equation}

Now,
\begin{align}
&\int_\Gamma (H_0^-F)^2(\gamma)\mu(d\gamma) \notag\\
&\quad=\int_\Gamma \bigg(-\int_{\mathbb R^d}\gamma(dx) \exp[sE(x,\gamma\setminus x)](D_{x}^{-}F)(\gamma)\bigg)^2 \mu(d\gamma) \notag\\
&\quad=\int_\Gamma \mu(d\gamma)\bigg(-\int_{\mathbb R^d}\gamma(dx) \exp[sE(x,\gamma\setminus x)]
(e^{\langle \varphi,\gamma\rangle -\varphi(x)}-e^{\langle \varphi,\gamma\rangle})\bigg)^2 \notag\\
&\quad=\int_\Gamma \mu(d\gamma)\bigg(-\int_{\mathbb R^d}\gamma(dx) \exp[sE(x,\gamma\setminus x)]
e^{\langle \varphi,\gamma\rangle}(e^{-\varphi(x)}-1)\bigg)^2 \notag\\
&\quad=\int_\Gamma \mu(d\gamma)e^{\langle 2\varphi,\gamma\rangle}\int_{\mathbb R^d}\gamma(dx_1)\int_{\mathbb R^d}\gamma(dx_2)
\exp[sE(x_1,\gamma\setminus x_1)]\exp[sE(x_2,\gamma\setminus x_2)] \notag\\
&\qquad\times
(e^{-\varphi(x_1)}-1)(e^{-\varphi(x_2)}-1) \notag\\
&\quad=\int_\Gamma \mu(d\gamma)\int_{\mathbb R^d}z\,dx\,e^{-E(x,\gamma)}e^{\langle 2\varphi,\gamma\cup x \rangle}
\exp[2sE(x,\gamma)](e^{-\varphi(x)}-1)^2
 \notag\\
&\qquad\text{}+\int_\Gamma \mu(d\gamma)\int_{\mathbb R^d}z\,dx_1\int_{\mathbb R^d}z\,dx_2 \exp[-E(x_1,\gamma)-E(x_2,\gamma)
-\phi(x_1-x_2)]e^{\langle 2\varphi,\gamma\cup x_1 \cup x_2 \rangle} \notag\\
&\qquad\times \exp[sE(x_1,\gamma \cup x_2)]\exp[sE(x_2,\gamma \cup x_1)](e^{-\varphi(x_1)}-1)(e^{-\varphi(x_2)}-1) \notag\\
&\quad=\int_{\mathbb R^d}z\,dx(1-e^{\varphi(x)})^2\int_\Gamma \mu(d\gamma)\exp[(2s-1)E(x,\gamma)+
\langle 2\varphi,\gamma \rangle] \notag\\
&\qquad+\int_{\mathbb R^d}z\,dx_1\int_{\mathbb R^d}z\,dx_2\,e^{\varphi(x_1)}(1-e^{\varphi(x_1)})e^{\varphi(x_2)}
(1-e^{\varphi(x_2)}) \notag\\
&\qquad\times\exp[(2s-1)\phi(x_1-x_2)]\int_\Gamma \mu(d\gamma)\exp[(s-1)E(x_1,\gamma)
+(s-1)E(x_2,\gamma)+\langle 2\varphi,\gamma \rangle]\label{E}  
\end{align}
and 
\begin{align}
&\int_\Gamma (H_0^+F)^2(\gamma)\mu(d\gamma) \notag\\
&\quad=\int_{\mathbb R^d}z\ dx_1\,\int_{\mathbb R^d}z\,dx_2\,(e^{\varphi(x_1)}-1)
(e^{\varphi(x_2)}-1) \int_\Gamma \mu(d\gamma)\exp[(s-1)E(x_1,\gamma)\notag\\
&\qquad+(s-1)E(x_2,\gamma)+\langle 2\varphi,\gamma \rangle].\label{F}
\end{align}

Completely analogously to \eqref{E} we have
\begin{align}
&\int_\Gamma (H_\varepsilon^-F)^2(\gamma)\mu(d\gamma) \notag\\
&\quad=\int_{\mathbb R^d}z\,dx\int_{\mathbb R^d}dy_1\int_{\mathbb R^d}dy_2(e^{\varphi(x)}-1)^2 a_\varepsilon (x-y_1)
a_\varepsilon(x-y_2)\int_\Gamma \mu(d\gamma) \notag\\
&\qquad\exp[(2s-1)E(x,\gamma)+(s-1)E(y_1,\gamma)+(s-1)
E(y_2,\gamma)+\langle 2\varphi,\gamma \rangle] \notag\\
&\qquad+\int_{\mathbb R^d}z\,dx_1\int_{\mathbb R^d}z\,dx_2\int_{\mathbb R^d}dy_1\int_{\mathbb R^d}dy_2 \  
e^{\varphi(x_1)}(e^{\varphi(x_1)}-1)e^{\varphi(x_2)}(e^{\varphi(x_2)}-1) \notag\\
&\qquad\times a_\varepsilon (x_1-y_1)a_\varepsilon (x_2-y_2)\  
\exp[(2s-1)\phi(x_1-x_2)+(s-1)\phi(y_1-x_2)
\notag \\
&\qquad+(s-1)\phi(x_1-y_2)]\int_\Gamma \mu(d\gamma)\exp[(s-1)E(x_1,\gamma)+(s-1)E(x_2,\gamma) \notag\\
&\qquad+(s-1)E(y_1,\gamma)+(s-1)
E(y_2,\gamma)+\langle 2\varphi,\gamma \rangle].\label{H}
\end{align}
Let us make the change of variables $$y_1'=\varepsilon (y_1-x),\quad \quad y_2'=\varepsilon (y_2-x)$$ in the first integral, and
$$y_1'=\varepsilon (y_1-x_1),\quad \quad y_2'=\varepsilon (y_2-x_2)$$ in the second integral. Then omitting the primes in the 
notations of variables, we continue \eqref{H} as follows: 
\begin{align}   
&=\int_{\mathbb R^d}z\,dx\int_{\mathbb R^d}dy_1\int_{\mathbb R^d}dy_2(e^{\varphi(x)}-1)^2 a(y_1)a(y_2)
\int_\Gamma \mu(d\gamma)\exp[(2s-1)E(x,\gamma) \notag\\
&\quad+(s-1)E((y_1/\varepsilon)+x,\gamma)+(s-1)
E((y_2/\varepsilon)+x,\gamma)+\langle 2\varphi,\gamma \rangle] \notag\\
&\quad+\int_{\mathbb R^d}z\,dx_1\int_{\mathbb R^d}z\,dx_2\int_{\mathbb R^d}dy_1\int_{\mathbb R^d}dy_2 \  
e^{\varphi(x_1)}(e^{\varphi(x_1)}-1)e^{\varphi(x_2)}(e^{\varphi(x_2)}-1) \notag\\
&\quad\times a(y_1)a(y_2)\  \exp[(2s-1)\phi(x_1-x_2)+(s-1)\phi((y_1/\varepsilon)+x_1-x_2)
\notag \\
&\quad+(s-1)\phi((y_2/\varepsilon)+x_2-x_1)\int_\Gamma \mu(d\gamma)\exp[(s-1)E(x_1,\gamma)+(s-1)E(x_2,\gamma) \notag\\
&\quad+(s-1)E((y_1/\varepsilon)+x_1,\gamma)+(s-1)
E((y_2/\varepsilon)+x_2,\gamma)+\langle 2\varphi,\gamma \rangle].\label{G}
\end{align}

Next,
\begin{align}
&\int_\Gamma (H_\varepsilon^+F)^2(\gamma)\mu(d\gamma) \notag\\
&\quad=\int_\Gamma \mu(d\gamma)\bigg( -\int_{\mathbb R^d} \gamma(dx)\int_{\mathbb R^d}dy \ a_\varepsilon(x-y) \notag\\
&\qquad \times \exp[sE(x,\gamma \setminus x)-(1-s) E(y,\gamma\setminus x)](F(\gamma\setminus x \cup y)
-F(\gamma \setminus x))\bigg)^2 \notag\\
&\quad=\int_\Gamma \mu(d\gamma)\bigg( \int_{\mathbb R^d} \gamma(dx)\int_{\mathbb R^d}dy \ a_\varepsilon(x-y) \notag\\
&\qquad \times \exp[sE(x,\gamma \setminus x)-(1-s) E(y,\gamma\setminus x)]
(e^{\langle \varphi,\gamma\setminus x \rangle +\varphi(y)}-e^{\langle \varphi,\gamma\setminus x \rangle})\bigg)^2 \notag\\
&\quad=\int_\Gamma \mu(d\gamma)
\bigg( \int_{\mathbb R^d} \gamma(dx_1)\int_{\mathbb R^d} \gamma(dx_2)
e^{\langle \varphi,\gamma\setminus x_1 \rangle}e^{\langle \varphi,\gamma\setminus x_2 \rangle}
\int_{\mathbb R^d}dy_1\int_{\mathbb R^d}dy_2 \notag\\
&\qquad\times a_\varepsilon(x_1-y_1)a_\varepsilon(x_2-y_2)
\exp[sE(x_1,\gamma \setminus x_1)-(1-s) E(y_1,\gamma\setminus x_1)]\notag\\
&\qquad\times\exp[sE(x_2,\gamma \setminus x_2)-(1-s) E(y_2,\gamma\setminus x_2)](e^{\varphi(y_1)}-1)(e^{\varphi(y_2)}-1) \notag\\
&\quad=\int_\Gamma \mu(d\gamma)
\int_{\mathbb R^d} \gamma(dx)e^{\langle 2\varphi,\gamma\setminus x \rangle}
\int_{\mathbb R^d}dy_1\int_{\mathbb R^d}dy_2 \notag\\
&\qquad\times a_\varepsilon(x-y_1)a_\varepsilon(x-y_2)
\exp[sE(x,\gamma \setminus x)-(1-s) E(y_1,\gamma\setminus x)]\notag\\
&\qquad\times\exp[sE(x,\gamma \setminus x)-(1-s) E(y_2,\gamma\setminus x)](e^{\varphi(y_1)}-1)(e^{\varphi(y_2)}-1) \notag\\
&\qquad +\int_\Gamma \mu(d\gamma)
\int_{\mathbb R^d} \gamma(dx_1)\int_{\mathbb R^d} (\gamma \setminus x_1)(dx_2)
e^{\langle \varphi,\gamma\setminus x_1 \rangle}e^{\langle \varphi,\gamma\setminus x_2 \rangle}
\int_{\mathbb R^d}dy_1\int_{\mathbb R^d}dy_2 \notag\\
&\qquad\times a_\varepsilon(x_1-y_1)a_\varepsilon(x_2-y_2)
\exp[sE(x_1,\gamma \setminus x_1)-(1-s) E(y_1,\gamma\setminus x_1)]\notag\\
&\qquad\times\exp[sE(x_2,\gamma \setminus x_2)-(1-s) E(y_2,\gamma\setminus x_2)](e^{\varphi(y_1)}-1)(e^{\varphi(y_2)}-1) \notag\\
&\quad=\int_\Gamma \mu(d\gamma)\int_{\mathbb R^d}z\,dx \exp[-E(x,\gamma)]e^{\langle 2\varphi,\gamma \rangle}
\int_{\mathbb R^d}dy_1\int_{\mathbb R^d}dy_2 \notag\\
&\qquad\times\varepsilon^{2d}a(\varepsilon(x-y_1))a(\varepsilon(x-y_2)) 
\exp[sE(x,\gamma)-(1-s) E(y_1,\gamma)]\notag\\
&\qquad\times\exp[sE(x,\gamma)-(1-s) E(y_2,\gamma)](e^{\varphi(y_1)}-1)(e^{\varphi(y_2)}-1) \notag\\
&\qquad+\int_\Gamma \mu(d\gamma)\int_{\mathbb R^d}z\,dx_1\int_{\mathbb R^d}z\,dx_2\int_{\mathbb R^d}dy_1\int_{\mathbb R^d}dy_2 \notag\\
&\qquad\times\exp[-E(x_1,\gamma)-E(x_2,\gamma)-\phi(x_1-x_2)]e^{\langle \varphi,\gamma \cup x_2 \rangle}
e^{\langle \varphi,\gamma \cup x_1 \rangle}\notag \\
&\qquad\times a_\varepsilon(x_1-y_1)a_\varepsilon(x_2-y_2)
\exp[sE(x_1,\gamma \cup x_2)-(1-s) E(y_1,\gamma\cup x_2)]\notag\\
&\qquad\times \exp[sE(x_2,\gamma \cup x_1)-(1-s) E(y_2,\gamma \cup x_1)](e^{\varphi(y_1)}-1)(e^{\varphi(y_2)}-1) \notag\\
&\quad=:{\rm I}+{\rm II}. \label{I}
\end{align}

In the first integral in \eqref{I} let us make the change of variables 
$$y_1'=\varepsilon (y_1-x),\quad  y_2'=\varepsilon (y_2-x).$$
Then omitting the primes in the 
notations of variables, we continue {\rm I} as follows: 
\begin{align}
{\rm I}&=\int_{\mathbb R^d}dy_1 \int_{\mathbb R^d}dy_2\, a(y_1)a(y_2)
(e^{\varphi((y_1/\varepsilon)+x)}-1)
(e^{\varphi((y_2/\varepsilon)+x)}-1) \notag\\
&\quad\times\int_\Gamma \mu(d\gamma)\int_{\mathbb R^d}z\,dx \,
\exp[(2s-1)E(x,\gamma)+(s-1)E((y_1/\varepsilon)+x,\gamma) \notag\\
&\quad
 +(s-1)E((y_2/\varepsilon)+x,\gamma)
+\langle 2\varphi,\gamma \rangle].\notag 
\end{align}
Let us take $$x'=x+(y_1/\varepsilon),$$ then omitting the primes in the notations of variables, we get:
\begin{align}
{\rm I}&=\int_{\mathbb R^d}dy_1 \int_{\mathbb R^d}dy_2\, a(y_1)a(y_2)\int_{\mathbb R^d}z\,dx (e^{\varphi(x)}-1)
(e^{\varphi(x+((y_2-y_1)/\varepsilon))}-1)  a(y_1)a(y_2)\notag\\
&\quad\times\int_\Gamma \mu(d\gamma)\exp[(2s-1)E(x-(y_1/\varepsilon),\gamma) \notag\\
&\quad
+(s-1)E(x,\gamma)+(s-1)
E(x+((y_2-y_1)/\varepsilon),\gamma)+\langle 2\varphi,\gamma \rangle] .\label{J}
\end{align}

In the second integral in \eqref{I}, let us make the change of variables 
$$x_1'=\varepsilon (x_1-y_1),\quad \quad x_2'=\varepsilon (x_2-y_2).$$
Then omitting the primes, we have:
\begin{align}
{\rm II}&=\int_{\mathbb R^d}z\,dx_1\int_{\mathbb R^d}z\,dx_2\int_{\mathbb R^d}dy_1\int_{\mathbb R^d}dy_2\,   
e^{\varphi((x_1/\varepsilon)+y_1)}
e^{\varphi((x_2/\varepsilon)+y_2)}(e^{\varphi(y_1)}-1)(e^{\varphi(y_2)}-1) \notag\\
&\quad\times a(x_1)a(x_2)\  \exp[(2s-1)\phi(((x_1-x_2)/\varepsilon)+y_1-y_2)+(s-1)\phi((x_1/\varepsilon)+y_1-y_2)
\notag \\
&\quad+(s-1)\phi(y_1-y_2-(x_2/\varepsilon))]
\int_\Gamma \mu(d\gamma)\exp[(s-1)E(y_1,\gamma)+(s-1)E(y_2,\gamma) \notag\\
&\quad+(s-1)E((x_1/\varepsilon)+y_1,\gamma)+(s-1)
E((x_2/\varepsilon)+y_2,\gamma)+\langle 2\varphi,\gamma \rangle].\label{K}
\end{align}
Using \eqref{J} and \eqref{K}, we get 
\begin{align}
&\int_\Gamma (H_\varepsilon^+F)^2(\gamma)\mu(d\gamma) \notag\\
&=\int_{\mathbb R^d}z\,dx\int_{\mathbb R^d}dy_1\int_{\mathbb R^d}dy_2\,(e^{\varphi(x)}-1)
(e^{\varphi(x+((y_2-y_1)/\varepsilon)}-1)  a(y_1)a(y_2)\notag\\
&\quad\times
\int_\Gamma \mu(d\gamma)\exp[(2s-1)E(x-(y_1/\varepsilon),\gamma) \notag\\
&\quad+(s-1)E(x,\gamma)+(s-1)
E(x+((y_2-y_1)/\varepsilon),\gamma)+\langle 2\varphi,\gamma \rangle] \notag\\
&\quad+\int_{\mathbb R^d}z\,dx_1\int_{\mathbb R^d}z\,dx_2\int_{\mathbb R^d}dy_1\int_{\mathbb R^d}dy_2\,   
e^{\varphi((x_1/\varepsilon)+y_1)}
e^{\varphi((x_2/\varepsilon)+y_2)}(e^{\varphi(y_1)}-1)(e^{\varphi(y_2)}-1) \notag\\
&\quad\times a(x_1)a(x_2)\,  \exp[(2s-1)\phi(((x_1-x_2)/\varepsilon)+y_1-y_2)+(s-1)\phi((x_1/\varepsilon)+y_1-y_2)
\notag \\
&\quad+(s-1)\phi(y_1-y_2-(x_2/\varepsilon))
\int_\Gamma \mu(d\gamma)\exp[(s-1)E(y_1,\gamma)+(s-1)E(y_2,\gamma) \notag\\
&\quad+(s-1)E((x_1/\varepsilon)+y_1,\gamma)+(s-1)
E((x_2/\varepsilon)+y_2,\gamma)+\langle 2\varphi,\gamma \rangle].\label{L}
\end{align}
Completely analogously, we get
\begin{align}
&\int_\Gamma (H_0^-F)(\gamma)(H_\varepsilon^-F)(\gamma)\mu (d\gamma)\notag\\
&=\int_{\mathbb R^d}z\,dx\int_{\mathbb R^d}dy \,(e^{\varphi(x)}-1)^2 a(y)
\int_\Gamma \mu(d\gamma) \exp[(2s-1)E(x,\gamma)+(s-1)E((y/\varepsilon)+x,\gamma)
+\langle 2\varphi,\gamma \rangle]\notag\\
&\quad+\int_{\mathbb R^d}z\,dx_1\int_{\mathbb R^d}z\,dx_2\int_{\mathbb R^d}dy\,(e^{\varphi(x_1)}-1)(e^{\varphi(x_2)}-1)
e^{\varphi(x_1)}e^{\varphi(x_2)}a(y) \notag\\
&\quad\times \exp[(2s-1)\phi(x_1-x_2)+(s-1)\phi((y/\varepsilon)+x_1-x_2)]\notag\\
&\quad\times\int_\Gamma \mu(d\gamma)\exp[(s-1)E(x_1,\gamma)+(s-1)E(x_2,\gamma) 
+(s-1)E((y/\varepsilon)+x_1,\gamma)+\langle 2\varphi,\gamma \rangle],\label{Q}
\end{align}
and 
\begin{align}
&\int_\Gamma (H_0^+F)(\gamma)(H_\varepsilon^+F)(\gamma)\mu (d\gamma)\notag\\
&=\int_{\mathbb R^d}z\,dx_1\int_{\mathbb R^d}z\,dx_2 \int_{\mathbb R^d} dy \, a(x_2)(e^{\varphi(x_1)}-1)\notag\\
&\quad\times  (e^{\varphi(y)}-1)\int_\Gamma \mu(d\gamma)e^{\langle 2\varphi,\gamma \rangle}
\exp[(s-1)E((x_2/\varepsilon)+y,\gamma)+(s-1)E(y,\gamma)\notag\\ 
&\quad+ 
(s-1)E(x_1,\gamma)+(s-1)\phi((x_2/\varepsilon)+y-x_1)].\label{R}
\end{align}

Using the Ruelle bound and Lemma \ref{n}, we conclude that the integral over $\Gamma$ in the right hand side of equalities
\eqref{G}, \eqref{L}--\eqref{R}, are bounded by a constant, which is indepeandent of $\varepsilon$. Therefore, by the 
dominated convergence theorem, to find the limit of \eqref{G}, \eqref{L}--\eqref{R} as $\varepsilon \to 0$
it suffices to find the point-wise limit of the functions appearing before the integral over $\Gamma$, as well as the 
limit of the integrals over $\Gamma$ for fixed $x$ 
 ($x_1$ and $x_2$ respectively), $y_1$ and $y_2$.

To find the latter limits, we use Lemms \ref{T}. Then, using \eqref{2}, we see that \eqref{G} and \eqref{Q}
 converge to \eqref{E}, whereas \eqref{L} and \eqref{R} converge to \eqref{F}. Therefore, \eqref{D} and \eqref{S} hold.\quad $\square$


\end{document}